\documentclass[12pt,a4j]{article}
\usepackage{amssymb}

\begin{document}

\begin{center}
{\bf\large Efficiency and influence function of
estimators for ARCH models
}
\end{center}
\begin{center}
{\sc {\sc S. Ajay
Chandra}\footnote[2]{a.chandra@latrobe.edu.au}}
\end{center}
\begin{center}
{\it La Trobe University, Victoria, Australia}
\end{center}
\newcommand{\qed}{\hbox{\rule{6pt}{6pt}}}
\renewcommand{\abstractname}{{}}
 \begin{abstract}
\noindent{\bf Abstract.} This paper proposes a
closed-form optimal estimator based on the theory
of estimating functions for a class of linear
ARCH models. The estimating function (EF)
estimator has the advantage over the widely used
maximum likelihood (ML) and quasi-maximum
likelihood (QML) estimators that (i) it can be
easily implemented, (ii) it does not depend on a
distributional assumption for the innovation, and
(iii) it does not require the use of any
numerical optimization procedures or the choice
of initial values of the conditional variance
equation. In the case of normality,
the asymptotic distribution of the ML and QML
estimators naturally turn out to be identical
and, hence, coincides with ours. Moreover, a
robustness property of the EF estimator is
derived by means of influence function.
Simulation results show that the
efficiency benefits of our estimator relative to
the ML and QML estimators are substantial for
some ARCH innovation distributions.
\end{abstract}

\noindent\textbf{1. Introduction}\\\\
 Since the
seminal papers of Engle (1982) and Bollerslev
(1986), autoregressive conditional
heteroskedasticity (ARCH) and generalised ARCH
(GARCH) models have been proposed for modelling
time series with non-constant conditional
volatility.
Since then, these models have become perhaps the
most popular and extensively studied financial
econometric models (see e.g., Engle (1995);
Gouri\'{e}roux (1997); Mikosch (2003); Francq and
Zako\"{\i}an (2004)). The literature on the
subject is so vast that we will restrict
ourselves to directing the reader to fairly
comprehensive reviews by Bollerslev et al. (1992)
and Shephard (1996). An excellent survey of the
GARCH methodology in finance is also available
such as
Bauwens et al. (2006).\\
\indent ARCH model estimation can be achieved
using a variety of methods such as conditional
least squares (LS) estimation (Tj\o stheim
(1986)), maximum likelihood (ML) estimation under
the assumption of conditional normality,
quasi-maximum likelihood (QML) estimation (Weiss
(1986); Francq and Zako\"{\i}an (2004)),
generalized method of
moments (GMM) estimation (e.g., Rich et al.
(1991)). As is well-known, LS, QML and GMM
estimation methods yield inefficient and possibly
biased estimates relative to ML estimators when
the true innovation distribution is known (see
for example, Li and Turtle (2000)). However, the
possibility for misspecification of the
likelihood function
for ML and QML estimators motivates our
investigation of an alternative estimation method
for ARCH models.
\\
\indent The purpose of this paper is to propose
an estimator based on the estimating function
(EF) approach for ARCH models that improves
efficiency without any distributional assumptions
for the innovation. This EF estimator admits a
closed-form expression which is computationally
simple and compares favorably with the ML and QML
estimators. Moreover, the EF estimator naturally
turns out to have the same limiting distribution
as the ML and QML estimators and hence is also
fully efficient when the innovation distribution
is Gaussian. It is interesting to note that many
standard results in the estimation of ARCH models
based on conditional normality are
recoverable under the EF approach.\\
\indent The rest of the paper is organized as
follows. Section 2 describes
the conditional LS  and EF estimation procedures
 of ARCH models. In addition,
the asymptotic efficiency of the EF estimator
relative to the LS, ML and QML estimators is
discussed. In particular, the lower bound of the
asymptotic variance of the LS estimator is
formulated. In Section 3, a robustness of the EF
estimator is studied by means of influence
function. In Section 4, we perform an experiment
to examine the asymptotic behavior of EF, LS, ML
and QML estimators in terms of mean square errors
in a small and a large-sample of observations.
The study demonstrates that the efficiency
benefits of the EF estimator relative to the LS,
ML and QML estimators are substantial
for some ARCH innovation distributions.\\\\
\textbf{2. Estimating
function formulation and efficiency}\\\\
 In this
section, we describe the problem of estimation
for a class of ARCH($p$) models
characterized by
the equations
\begin{equation}
X_t=
\sigma_{t}(\theta_0)\varepsilon_t,\quad
\sigma_{t}^2 (\theta_0)=
\omega_0+\sum\limits_{j=1}^p \alpha_{0j}
X_{t-j}^2,\quad t=1,\ldots, n,
\end{equation}
where {$\{\varepsilon_t\}$ is a sequence of
independent, identically distributed random
variables such that $E\varepsilon_t=0 $,
$E\varepsilon_t^2=1$,
$\theta_0=(\omega_0,\alpha_{01}
\ldots,\alpha_{0p})^T
$ is an unknown vector of true parameters
satisfying $\omega_0>0$, $\alpha_{0j}\ge 0$,
$j=1,\ldots,p-1$, and $\varepsilon_t$ is
independent of $X_s,s<t$.
Henceforth, it is tacitly assumed $\alpha_{0p}>0$
so that the model is of order $p$. We also assume
that model (1) is stationary and ergodic. When
$p=1$, Nelson (1990) showed that a sufficient
condition for the stationarity is $E(\log(
\alpha_{01} \varepsilon_t^2))<0$. For a general
ARCH($p$) model, Bougerol and Picard (1992)
showed that it has a unique non-anticipative
strictly stationary solution.

We now turn to describe the conditional least
squares estimation of model (1). Write
$Y_t=(1,X_t^2,\ldots,X_{t-p+1}^2)^T$ and
$\eta_t=(\varepsilon_{t}^2-1)\sigma_t^2(\theta_0)$.
Then the standard linear autoregressive
representation
is given by
\begin{equation} X_t^2=\theta_0^T
Y_{t-1}+\eta_t.
\end{equation}
Suppose that an observed stretch
$X_1^2,\ldots,X_n^2$ is available. The vector of
parameters is $\theta=(\omega, \alpha_1
\ldots,\alpha_{p})^T$ which belongs to a compact
parameter space $\Theta \subset (0,\infty)\times
[0,\infty)^p$, and $\theta_0\in \Theta$. Let
$$Q_n(\theta)=\sum_{t=1}^n
(X_t^2-E(X_t^2|\mathcal{F}_{t-1}))^2=\sum_{t=1}^n
(X_t^2-\theta^T Y_{t-1})^2$$ be the penalty
function, where $\mathcal{F}_t=\sigma\{X_s^2,t\le
s\}$. Then from the linear regression theory, we
can define the conditional least squares (LS)
estimator of $\theta$ by
\begin{equation}
 \hat{\theta}_n^{(LS)}=\arg \min_{\theta \in
 \Theta}
 Q_n(\theta)=(Y^TY)^{-1}YX,
\end{equation}
where $Y$ is the matrix of order $n\times (1+p)$
with $t$th row $Y_{t-1}$ and
$X=(X_1^2,\ldots,X_n^2)^T$.
Note that (3) does not take into account the
nuisance parameter $\theta_0$ associated with
variance and, hence, serves only as an initial
estimator. The asymptotic validity of (3) can be
easily established using an appropriate central
limit theorem, and as expected, its efficiency is
smaller than that of the QML estimator.
 \\
\indent To describe the limiting distribution of
(3), we impose an additional condition on
$\theta_0$, and the moment of $\varepsilon_t$.
Recall that we have assumed that the process is
stationary and ergodic. Let
\[
\mathcal{A}_{0t}=\left(
\begin{array}{cccc}
\alpha_{01} \varepsilon_t^2&\cdots&
\alpha_{0p-1}\varepsilon_t^2&
\alpha_{0p}\varepsilon_t^2\\
1&\cdots&0&0\\
\vdots&\ddots&\vdots&\vdots\\
0&\cdots&1&0\\
\alpha_{01} &\cdots& & \alpha_{0p}
\end{array}
\right).
\]
Introduce the notation $\mathcal{A}_{0t}^{\otimes
s} =\mathcal{A}_{0t}\overbrace{\otimes \cdots
\otimes}^{s \;{\rm times}} \mathcal{A}_{0t}$,
$\Sigma_s=E(\mathcal{A}_{0t}^{\otimes s})$, where
$\otimes$ denotes the tensor product.
\\\\
\textbf{Assumption 1}
\begin{itemize}
\item[] $E|\varepsilon_t|^8<\infty$ and
$\|\Sigma_{4}\|<1$,
\end{itemize}
where $\|\cdot\|$ is the spectral matrix norm.
In the case when $p=1$, and $\{\varepsilon_t\}$
is Gaussian, it is seen that $\|\Sigma_4\|<1$
implies $\alpha_{01}<105^{-\frac14}\approx 0.3$.
The following theorem establishes the asymptotic
distribution of (3).\\\\
\textbf{Lemma 1}. {\it Suppose that the
assumptions of model (1) and Assumption 1 hold.
Then
\[
\sqrt{n}
(\hat{\theta}_n^{(LS)}-\theta_0)\stackrel{d}{\to}
\mathcal{N}(0,Var(\varepsilon_t^2)\mathcal{U}^{-1}
\mathcal{R}(\theta_0)\mathcal{U}^{-1}) \quad
\mbox{as $n\to\infty$},
\]
where the matrices $\mathcal
{U}=E(Y_{t-1}Y_{t-1}^T)$ and
$\mathcal{R}(\theta_0)= E(\sigma_t^4(\theta_0)
 Y_{t-1}
Y_{t-1}^T)$ are positive definite with typically
bounded
elements.}\\\\
\textbf{Remark 1.} Note that Assumption 1 ensures
that $\mathcal {U}$ and $\mathcal{R}(\theta_0)$
are all finite. When the errors are standard
normal, necessary and sufficient condition for
the existence of higher moments of $X_t$ in terms
of the parameter $\theta_0$ is given by Engle
(1982, Theorems 1 and 2).\\\\
\noindent\textbf{Remark 2.} As an illustration,
we verify $\mathcal{R}(\theta_0)$ is positive
definite. Indeed, it is nonnegative definite,
i.e., $c^T \mathcal{R}(\theta_0)c=E(
\sigma_t^2(\theta_0)c^T
 Y_{t-1})^2 \ge 0$ for any given vector $c=
 (c_0,\ldots,c_p)^T\in R^{p+1}$.
 Moreover, if we suppose that $\mathcal{R}(\theta_0)$ is
 not positive definite, then there exists a
 vector $(c_0,\ldots,c_{j_0})$ with
 $c_{j_0}\not=0$ ($j_0 \le p$) such that $c_0+
 c_1 X_{s-1}^2+\cdots+ c_{j_0}X_{s-j_0}^2=0$
 a.e. This implies $\sigma_t^2(\theta_0)>0$ a.e.,
 because of $\omega_0>0$. In this case, we can
 write $X_{s-j_0}^2=-\gamma_0-\gamma_1X_{s-1}^2-\cdots-
 \gamma_{j_0-1}X_{s-j_0+1}^2$, where
 $\gamma_k=c_k/c_{j_0}$. Hence,
 substituting this into the last term of
 $\sigma_s^2(\theta_0)$ in (1) with $s-j_0=t-p$
 entails an
 ARCH($p-1$) representation,
 leading to a contradiction.\\\\
\indent The conditional least squares estimator
$\hat{\theta}_n^{(LS)}$ typically possesses the
properties that it admits a closed-form
expression, which is computationally easy.
However, $\hat{\theta}_n^{(LS)}$ in general is
not asymptotically efficient. Thus we next
discuss an asymptotically efficient estimator
proposed by Godambe (1985) which has the
following desirable properties: (i) it has an
explicit form which is computationally easy (ii)
it compares favourably with the ML and QML
estimators.

\indent Let $X^{(n)}=(X_1,\ldots,X_n)^T$ be a
vector of random variables forming a stochastic
process. The distribution family $\mathbb{F}$ of
$X^{(n)}$ is specified by an unknown parameter
vector $\theta=(\theta_1,\ldots,\theta_p)^T$ and
$\theta=\theta(F)$, $F \in {\mathbb F}$ be a real
parameter vector.
An estimating function $g( X^{(n)},\theta(F))$
satisfying certain regularity conditions, is
called a regular unbiased estimating function if
\begin{equation}
E[g(X^{(n)},\theta(F))]=0,\quad F\in{\mathbb F}.
\end{equation}
For a given set of unbiased estimating functions
$g$ belonging to the class $\mathcal{G}$, the
estimating function $g^*\in \mathcal{G}$ is said
to be optimal for $\theta$ if
\[
{\rm E}[g^2(X^{(n)},\theta(F))]/\{{\rm E}[
\partial g(X^{(n)},\theta(F))/\partial
\theta |_{\theta=\theta(F)}]\}^2
\]
is minimized for all $F\in{\mathbb F}$ at
$g=g^*$.\\
\indent Let $\mathbb{L}$ be a class of linear
combinations of unbiased estimating functions
of the
form
\[
g=\sum_{t=1}^n a_{t-1}h_t ,
\]
where the weights $a_{t-1}$ are any function of
$X_1, \ldots,X_{t-1}$ and $\theta$, and $h_t$ is
a function of $X_1,\ldots,X_t$ and $\theta$
satisfying
$E_F(h_t|{\cal B}_{t-1})=0$, where
$\mathcal{B}_t=\sigma\{X_s,t\le s\}$. Moreover,
for all $F\in \mathbb{F}$, the $h_t$'s are
mutually
orthogonal.\\
\indent An obvious example of $h_t$ is
\begin{equation}
h_t=X_t-E(X_t|{\mathcal B}_{t-1}),
\end{equation}
which is the residual between $X_t$ and its best
predictor $E(X_t|{\mathcal B}_{t-1})$.
 We assume that $h_t$ and $a_{t-1}$ are
differentiable with respect to $\theta$ for $1\le
t \le n$. These considerations motivate the
following result, which is due to Godambe (1985).\\

\noindent\textbf{Lemma 2}. {\it In the class
$\mathcal{G}$ of estimating functions $g$, the
optimal estimating function is given by
\[
g^*  =\sum_{t=1}^n a_{t-1}^* h_t,
\]
where
 $a_{t-1}^*=E[(\partial h_t/\partial \theta)|
 \mathcal{B}_{t-1}]/E[h_t^2|\mathcal{B}_{t-1}]$.}\\\\
 \indent By virtue of Lemmas 1 and 2, we are now in a
 position to state our main result. For this purpose,
 we need the following notation. In view of (5),
 we can set
 \begin{equation}
   h_t=X_t^2-E(X_t^2|\mathcal{F}_{t-1})
   =X_t^2-\sigma_t^2(\theta_0).
\end{equation}
In general, the choice of an estimating function
can be viewed in a manner analogous to the
selection of moment conditions in the generalized
method of moments (see Hansen (1982)). Now based
on (6), we have
\[
E\biggl(\frac{\partial h_t}{\partial
\theta}\biggl|\mathcal{F}_{t-1}\biggl)=-
\frac{\partial \sigma_t^2(\theta_0)}{\partial
\theta}\quad \mbox{and}\quad
E(h_t^2|\mathcal{F}_{t-1})=
E(X_t^4|\mathcal{F}_{t-1})-\sigma_t^4(\theta_0).
\]
Then by virtue of (6) and Lemma 2, the optimal
estimating
 function is
\begin{equation}
g^*=-\sum_{t=1}^n \frac{\frac{\partial
\sigma_t^2(\theta_0)}{\partial \theta}(X_t^2-\sigma_t^2(\theta_0))}
{E(X_t^4|\mathcal{F}_{t-1})-\sigma_t^4(\theta_0)}.
\end{equation}
It should be pointed out
that (7) is based on the finite sample, and it
does not depend on any distributional assumptions
for $X_t^2$ conditional on $\mathcal{F}_{t-1}$.
Noting that
$E(h_t^2|\mathcal{F}_{t-1})=Var(\varepsilon_t^2)
\sigma_t^4(\theta_0)$ and using (3),
it follows that the solution to $g^*=0$ in (7) is
the estimating function (EF) estimator given by
\begin{equation}
\hat{\theta}_n^{(EF)}=\left(\sum_{t=1}^n
\frac{Y_{t-1}
Y_{t-1}^T}{\sigma_t^4(\hat{\theta}_n^{(LS)})}\right)^{-1}
\left(\sum_{t=1}^n
\frac{Y_{t-1}X_t^2}{\sigma_t^4(\hat{\theta}_n^{(LS)})}\right).
\end{equation}
Here, it is assumed that
$\sqrt{n}(\hat{\theta}_n^{(LS)}-\theta_0)=O_p(1)$.
We now impose the following additional regularity
conditions. Recall the matrix $\mathcal{A}_{0t}$
and write it as
$\mathcal{A}_0=\{\mathcal{A}_{0t}\}$. In the
notation of Bougerol and Picard (1992),
the top Lyapunov exponent is defined by
$$\gamma(\mathcal{A}_0)\equiv \inf _{t\ge
1}\frac{1}{t} E(\log
\|\mathcal{A}_{01}\mathcal{A}_{02}\cdots
\mathcal{A}_{0t} \|)=\lim
_{t\to\infty}\frac{1}{t}\log
\|\mathcal{A}_{01}\mathcal{A}_{02}\cdots
\mathcal{A}_{0t} \|$$ under the assumption that
$E(\log^+\|\mathcal{A}_{01}\|)\le
E\|\mathcal{A}_{01}\|<\infty$.  \\\\
 \textbf{Assumption 2}
\begin{itemize}
\item[(i)] $\theta_{0}\in \tilde{\Theta}$, where
$\tilde{\Theta}$ denotes the interior of the
compact parameter space $\Theta$. \item[(ii)]
$\gamma(\mathcal{A}_{0})<0$.
 \item[(iii)] $\varepsilon_{t}^2$ has a non-degenerate
 distribution with $E \varepsilon_{t}^2=1$.
 \item[(iv)]
$E\varepsilon_{t}^4<\infty$. \item [(v)]
$\Gamma(\theta_0)=E(Y_{t-1}
Y_{t-1}^T/\sigma_t^4(\theta_0))$ is finite.
\end{itemize}    Hence, we have the following theorem,
which is the main result of the paper. The proofs
for Lemma 1 and Theorem 1 are given in Section
5.\\\\
 \noindent\textbf{Theorem 1.} {\it
 Suppose that the
assumptions of model (1) and Assumption 2 hold.
Then
as $n\to \infty$,
\[
\sqrt{n}( \hat{\theta}_n^{(EF)}-\theta_0)
\stackrel{d}{\rightarrow}\mathcal{
N}(0,Var(\varepsilon_t^2)\Gamma^{-1}(\theta_0)).
\]}

\noindent\textbf{Remark 3.} Under the assumption
of conditional normality, we have
$E(X_t^4|\mathcal{B}_{t-1})=3\sigma_t^4(\theta_0)$,
and in analogy with Engle (1982) the expression
(7) reduces to
\begin{equation}
g^*=-\sum_{t=1}^n \frac{1}{2\sigma_t^2(\theta_0)}
\frac{\partial \sigma_t^2(\theta_0)}{\partial
\theta}\biggl(\frac{X_t^2}{h_t}-1\biggl).
\end{equation}
Comparing (9) with the first-order condition of
Engle (1982), we observe that they are equivalent
up to a sign change. As is well known, under the
additional assumption of normality, the ML
estimator of $\theta$ has the following
asymptotic distribution,\[ \sqrt{n}(
\hat{\theta}_n^{(ML)}-\theta_0)
\stackrel{d}{\rightarrow}\mathcal{
N}(0,2\Gamma^{-1}(\theta_0)),\quad \mbox{as $n\to
\infty$}.
\]
Hence, we conclude that the theory of estimating
functions and ML method for the estimators of the
ARCH model yield essentially the same asymptotic
distribution.\\\\
\noindent\textbf{Remark 4.} As shown by Francq
and Zako\"{\i}an (2004) for the ARCH model, the
QML estimator of $\theta$ is obtained by
maximising the normal log-likelihood function
although the true probability density function is
non-normal. Under the conditions of model (1) and
Assumption 2, they showed that the QML estimator
is asymptotically normal:
\[
\sqrt{n}( \hat{\theta}_n^{(QML)}-\theta_0)
\stackrel{d}{\rightarrow}\mathcal{
N}(0,Var(\varepsilon_t^2)\Gamma^{-1}(\theta_0)),\quad
\mbox{as $n\to \infty$}.
\]
It is interesting to note that, if the true
probability density function is the normal
distribution, the asymptotic distribution of the
ML and QML estimators is identical and coincides
with ours.\\\\
\noindent\textbf{Remark 5.}
Consider the GARCH($p,q$) model
\begin{eqnarray}
X_t=
\sigma_{t}(\vartheta_0)\varepsilon_t,\quad
\sigma_{t}^2 (\vartheta_0)=
\omega_0+\sum\limits_{i=1}^p \alpha_{0i}
X_{t-i}^2+\sum\limits_{j=1}^q \beta_{0j}
\sigma_{t-j}^2 (\vartheta_0),\quad t=1,\ldots,n
\end{eqnarray}
where {$\{\varepsilon_t\}$ is a sequence of
independent, identically distributed random
variables such that $E\varepsilon_t=0$,
$E\varepsilon_t^2=1$,
$\vartheta_0=(\omega_0,\alpha_{01}
\ldots,\alpha_{0p},\beta_{01},\ldots,\beta_{0q})^T\in
\Theta \subset (0,\infty)\times [0,\infty)^{p+q}$
is an unknown vector of true parameters
satisfying $\omega_0>0$, $\alpha_{0i}\ge 0$,
$i=1,\ldots,p$, $\beta_{0j}\ge0$,
 $j=1,\ldots,q$ and $\varepsilon_t$ is
independent of $X_s,s<t$. Necessary and
sufficient conditions under which the
 GARCH($p,q$) equations have a unique, strictly
 stationary, and non-anticipative solution
 were found by Nelson (1990)
 for $p=1$ and $q=1$, and by Bougerol and Picard (1992a, b) for
 arbitrary $p\ge 1$ and $q\ge 1$.  \\
\indent Write
$\xi_t=(\varepsilon_t^2-1)\sigma_t^2(\vartheta_0)$.
Then by analogy with (2), it follows that (10)
can be
 represented as an ARMA($p^*,q$) model:
\begin{equation}
X_t^2= \omega_0 +\sum_{i=1}^{p^*}\phi_{0i}
X_{t-i}^2 +\xi_t-\sum_{j=1}^q\beta_{0j}\xi_{t-j},
\end{equation}
where $p^*=\max\{p,q\}$
and $\phi_{0i}=\alpha_{0i}+\beta_{0i}\ge 0$,
$i=1,\ldots,p^*$.
 We have further defined $\alpha_{0i}=0$ for $i>p$ and
 $\beta_{0j}=0$ for $j>q$.
 Henceforth, it is assumed that $X_t^2$ is
 covariance-stationary
  provided that $\xi_t$ has finite variance
  and that the roots
  of $1-\phi_{01}z
  -\cdots-
  \phi_{0p^*}z^{p^*}=0$ are outside the
  unit circle. Given the nonnegativity
  restriction, this means that $X_t^2$ is
  covariance-stationary if $\phi_{01}+\ldots+\phi_{0p^*}
  <1$.
   \\
\indent
Suppose that an observed stretch
$X_1^2,\ldots,X_n^2$ is available from
$\{X_t^2\}$.
Let $R(l)$ denote the autocovariance function of
lag $l$,
\[
R(l)=E{(X_t^2-\mu)(X_{t+l}^2-\mu)}, \quad
l=0,\pm1,\ldots,
\]
where $\mu=
E(X_t^2)=\omega_0/(1-\phi_{01}-\cdots-\phi_{0p^*})
$, with the corresponding estimator
\[
\hat{R}(l)= \frac{1}{n}
\sum_{t=1}^{n-|l|}(X_t^2-\hat{\mu})(X_{t+|l|}^2-
\hat{\mu}),\quad
|l|<n,
\]
where $\hat{\mu}=\sum_{t=1}^n X_t^2/n$. Here the
initial values of
$X_0^2=\cdots=X_{1-p^*}^2=\xi_0=
 \cdots=\xi_{1-q}=0$ have negligible effect on parameter
 estimates when the sample size is large. The
expression
\[
I_n(\lambda)= \frac{1}{2 \pi
n}\left|\sum_{t=1}^n(X_t^2-\hat{\mu})e^{it
\lambda}\right|^2 = \frac{1}{2
\pi}\sum_{l=-n+1}^{n-1}\hat{R}(l)e^{-il \lambda}
\]
is called the periodogram of the partial
realization of $\{X_t^2\}$. The vector of
parameters is $\varphi=(\omega, \phi_1
\ldots,\phi_{p^*},\beta_1,\ldots,\beta_q)^T\in
\Theta$. By the stationarity and ergodicity, set
 $\sigma_\xi^2=E\xi_t^2>0$ and write
 $\rho=(\varphi^T,\sigma_\xi^2)^T$. Then the
 spectral density
 of $\{X_t^2\}$ is
\[
f_\rho(\lambda)= \frac{\sigma_\xi^2}{2
\pi}\biggl|1-\sum\limits_{j=1}^q\beta_j
e^{ij\lambda}\biggr|^2\biggl|1-\sum\limits_{j=1}^{p^*}\phi_j
e^{ij\lambda}\biggr|^{-2}.
\]
In order to estimate $\rho$,
Hosoya and Taniguchi (1982) proposed to minimise
\[
D(f_\rho(\lambda), I_n)= \int_{-\pi}^\pi
\biggl\{\log f_\rho(\lambda)+\frac{
f_\rho(\lambda)}{ I_n(\lambda)}\biggl\} d\lambda
\]
with respect to $\rho$.
Let $\hat{\rho}^{(QML)}=
(\hat{\omega}^{(QML)},\hat{\phi}_1^{(QML)},
\ldots,\hat{\phi}_{p^*}^{(QML)},
\hat{\beta}_1^{(QML)},\ldots,\hat{\beta}_q^{(QML)})^T$
be a quasi-Gaussian
 maximum likelihood estimator of $\rho$ which
 minimizes $D(f_\rho(\lambda), I_n)$. Under
 certain conditions, they showed that
\begin{eqnarray*}
\sqrt{n}(\hat{\rho}^{(QML)}-\rho)\stackrel{d}
{\longrightarrow}{\cal
N}\biggl(0,4\pi\biggl(\int_{-\pi}^\pi
\frac{\partial}{\partial \rho}\log
f_\rho(\lambda)\frac{\partial}{\partial
\rho^T}\log
f_\rho(\lambda)d\lambda\biggr)^{-1}\biggr).
\end{eqnarray*}
Note that (11) does not take into account the
nuisance parameter $\vartheta_0$ associated with
variance. Hence, $\hat{\rho}_n^{(QML)}$ indeed
serves only as an initial estimator. \\
\indent The vector of parameters is
$\vartheta=(\omega, \alpha_1
\ldots,\alpha_p,\beta_1,\ldots,\beta_q)^T\in
\Theta$. To apply Godambe's estimating function
method to (11), let
$$\hat{\sigma}_t^2\equiv
\sigma_t^2(\hat{\vartheta}^{(QML)})=\hat{\omega}^{(QML)}
+\sum_{i=1}^p
\hat{\alpha}_i^{(QML)}X_{t-i}^2+\sum_{j=1}^q
\hat{\beta}_j^{(QML)}\sigma_{t-j}^2
(\hat{\vartheta}^{(QML)}),
$$
with
$\hat{\sigma}_0^2=\hat{\sigma}_{-1}^2=\cdots=
\hat{\sigma}_{-q+1}^2=0$. Then we can construct
$\hat{\sigma}_t^2$, $t=1,\ldots,n$, iteratively.
Once this is done, we can find an estimator
$\hat{\hat{\rho}}_n^{(EF)}=\hat{\rho}_n^{(EF)}
(\hat{\rho}_n^{(QML)})$ of $\rho$ by means of
 Godambe's method. \\

\indent To gain a further insight into Theorem 1,
it is interesting to compare the asymptotic
variances of $\hat{\theta}_n^{(LS)}$ and
$\hat{\theta}_n^{(EF)}$ in terms of their
efficiency. This motivates to state the following
result, whose proof is given in
Section 4.\\\\
\textbf{Theorem 2.} {\it Under the conditions of
Lemmas 1 and 2, the asymptotic variance of
$\hat{\theta}_n^{(LS)}$ typically satisfies the
inequality that
$Var(\varepsilon_t^2)\mathcal{U}^{-1}
\mathcal{R}(\theta_0)\mathcal{U}^{-1} \ge
Var(\varepsilon_t^2)\Gamma^{-1}(\theta_0)$.}
\\ \\
       \noindent\textbf{3. Influence function}\\\\
    \noindent An influence function is
    a statistical tool which provides
    rich qualitative information of how an estimator
    responds to a small
    amount of contamination at any point. In the
    following we introduce a robustness measure
    for the estimator given by (8).
    This robust measure by means of influence function will
    show how the initial estimator
    $\hat{\theta}_n^{(LS)}$ affects
    $\hat{\theta}_n^{(EF)}$.\\
\indent Let us first study a robustness property
of $\hat{\theta}_n^{(LS)}$. Write
$S_t=(X_t^2,\ldots,X_{t-p+1}^2)^T$ and
$Z_{S,t}=(1,S_t^T)^T$. Then we have
$\hat{\theta}_n^{(LS)}=
\hat{\mathcal{U}}_S^{-1}\hat{\gamma}_S$, where
\[
\hat{\gamma}_S=\frac1n\sum_{t=1}^n S^{(1)}_t
Z_{S,t-1} \quad
\mbox{and}\quad\hat{\mathcal{U}}_S=\frac1n\sum_{t=1}^n
Z_{S,t-1} Z_{S,t-1}^T.
\]
Here, $S_t^{(1)}$ is the first component of
$S_t$. We can now define the corresponding least
squares functional as
$$T_{S}\equiv \mathcal{U}_{S}^{-1}\gamma_{S},$$
where
\[
\gamma_S =E(S_t^{(1)} Z_{S,t-1})\quad
\mbox{and}\quad
\mathcal{U}_S=E(Z_{S,t-1}Z_{S,t-1}^T),
\]
\indent As a measure of robustness property of
$\hat{\theta}_n^{(LS)}$, we consider the
following contaminated process
\[
S_{\delta,t}=(1-\delta)S_t+\delta U_t\equiv
S_t+\delta V_t,
\]
where $\delta\in(0,1)$. For
$S_\delta=\{S_{\delta,t}\}$, we can introduce an
influence function
\[
T_{S}'=\lim_{\delta\searrow0}\frac{T_{S_\delta}-T_S}
{\delta}.
\]
Noting the formula for differentiation of the
inverse of a matrix $d A^{-1}=-A^{-1}(d
A)A^{-1}$, we obtain
\begin{equation}
\frac{d}{d
\delta}\mathcal{U}_{S_\delta}^{-1}\biggl|_{\delta=0}
=-\mathcal{U}_S^{-1}
(\Delta+\Delta^T)\mathcal{U}_S^{-1},
\end{equation}
where $\Delta=E( \tilde{V}_{t-1}Z_{S,t-1}^T)$
with $\tilde{V}_t=(0,V_t^T)^T$.
Also,
\begin{equation}
\frac{d}{d
\delta}\gamma_{S_\delta}\biggl|_{\delta=0}
=E(V_t^{(1)}Z_{S,t-1})+E(S_t^{(1)}\tilde{V}_{t-1})\equiv\gamma_{S}',
\end{equation}
where $V_t^{(1)}$ is the first component of
$V_t$. Hence,
\[
T_S'=\mathcal{U}_S^{-1}[\gamma_S'-
(\Delta+\Delta^T)T_S].
\]
The quantity $T_S'$ will reveal how outliers in
the dependent and independent variables may
combine to affect $\hat{\theta}_n^{(LS)}$. \\
\indent In the above notation, we can similarly
derive an influence function of
$\hat{\theta}_n^{(EF)}=
\hat{\mathcal{U}}_{S,w}^{-1}\hat{\gamma}_{S,w}$,
where
\[
\hat{\gamma}_{S,w}=\frac1n\sum_{t=1}^n S^{(1)}_t
\hat{w}_{S,t-1}
 Z_{S,t-1} \quad
\mbox{and}\quad\hat{\mathcal{U}}_{S,w}=\frac1n\sum_{t=1}^n
\hat{w}_{S,t-1} Z_{S,t-1} Z_{S,t-1}^T
\]
with $$\hat{w}_{S,t}=[w((\hat{\theta}_n^{(LS)})^T
Z_{S,t} Z_{S,t}^T\hat{\theta}_n^{(LS)})]^{-1}.$$
Here note that
$\hat{w}_{S,t-1}=\sigma_t^{-4}(\hat{\theta}_n^{(LS)})$.
Since $\hat{\gamma}_{S,w}$ and
$\hat{\mathcal{U}}_{S,w}$ are the respective
sample versions of
 \[
\gamma_{S,w} =E[S_t^{(1)}w_{S,t-1}
Z_{S,t-1}]\quad \mbox{and}\quad
\mathcal{U}_{S,w}=E[w_{S,t-1}
Z_{S,t-1}Z_{S,t-1}^T],
\]
 with
 $$w_{S,t}
 =[w((T_S')^T
Z_{S,t} Z_{S,t}^TT_S')]^{-1},$$ the functional
analogue of $T_S$ is $$T_{S,w}=
 \mathcal{U}_{S,w}^{-1}\gamma_{S,w},$$  Write
$$Q_{t}=
\tilde{V}_tZ_{S,t}^T\quad \mbox{and}\quad
L_{w,t}= S_{t+1}^{(1)}w_{S,t}^2
 (T_S')^T Q_{t}T_S'.$$ Then by analogy
with (12),
 \[
\frac{d}{d
\delta}\mathcal{U}_{S_\delta,w}^{-1}\biggl|_{\delta=0}
=-\mathcal{U}_{S,w}^{-1}
[\tilde{\Delta}_w+\tilde{\Delta}_w^T-(\Phi_w+\Phi_w^T)
] \mathcal{U}_{S,w}^{-1},
\]
where
$$\tilde{\Delta}_w=E(w_{S,t-1}\tilde{V}_{t-1}Z_{S,t-1}^T
)\quad \mbox{and}\quad \Phi_w=E(L_{w,t-1}
Z_{S,t-1} Z_{S,t-1}^T )$$ and with (13),
 \begin{eqnarray*}
\frac{d}{d
\delta}\gamma_{S_\delta,w}\biggl|_{\delta=0}
&=&E(V_t^{(1)} w_{S,t-1}Z_{S,t-1})+E(S_t^{(1)}w_{S,t-1}
\tilde{V}_{t-1})\\
&&- E(L_{w,t-1} Z_{S,t-1}  )- E(L_{w,t-1}^T
Z_{S,t-1})\\
 &\equiv&\gamma_{S,w}'.
\end{eqnarray*}
Hence
 \[
T_{S,w}'=\mathcal{U}_{S,w}^{-1}\{\gamma_{S,w}'-
[\tilde{\Delta}_w+\tilde{\Delta}_w^T-(\Phi_w+\Phi_w^T)]T_{S,w}\}.
\]
 This expression will facilitate the fundamental
 description of sensitiveness or insensitiveness
 of $\hat{\theta}_n^{(EF)}$.\\\\
\noindent\textbf{4. Simulations}\\\\
 \noindent A finite sample experiment is performed
 to assess
 the asymptotic efficiency of the EF estimator
 given by (8) relative to LS, ML and QML
 estimators for a small and a large sample of
 observations.\\
\indent To facilitate meaningful comparisons, we
 generate an ARCH(1) of length $n=50$ or $n=500$
 for values of $\theta_0=(\omega_0,\alpha_{01})
 =(1,0.1)$ or
$(1,0.3)$
  based on 5000 replications. Without loss of efficiency,
  we assume $\omega_0=1$ and estimate
  $\alpha_{01}$ using LS, ML, QML and
  EF methods.
For this purpose,
we consider two error distributions$-$a normal
distribution and a Student$-$$t_v$ distribution
with $v=5$ degrees of freedom. In both cases, the
  unconditional mean and variance are 0 and 1,
  respectively. \\
\indent Tables 1 and 2 report the results in
terms of bias, variance and mean square error
  (MSE) for $\hat{\alpha}_n^{(LS)}$,
  $\hat{\alpha}_n^{(ML)}$,
  $\hat{\alpha}_n^{(QML)}$ and $\hat{\alpha}_n^{(EF)}$
  of $\alpha_{01}$. A closer examination of the MSE
  values in the tables reveals some interesting
  features. We first note that the values are
  intrinsically stable with respect to the choice
  of parameters and sample sizes. In every case,
  we observe that the finite-sample MSE of the EF
  estimator is desirable relative to other
  alternatives such as LS, ML and QML methods.
  The efficiency of $\hat{\alpha}_n^{(EF)}$
  increases against its counterpart as the sample
  size increases or $\alpha_{01}$ decreases. In the
  case of normality, the MSE
  results for the ML, QML and EF estimators are
  approximately the same. To this end, note that the
  result of
  Tables 1 and 2 also holds true for normalized
  error distributions such as double exponential,
  logistic and gamma having
  zero mean and unit variance. \\
  \indent Our simulation results highlight the
  benefits of using the EF formulation for
  modelling data drawn from non-normal
  conditional distributions. This approach
  naturally takes advantage of departures from
  normality to improve the efficiency of
  estimators given a finite sample of data. In
  comparison with asymptotically based methods,
  the focus on the finite sample in the EF
  formulation is important. We observe that
  efficiency gains from the EF approach are
  substantial. Hence, the EF formulation is
  potentially useful in cases with serious
  departures from normality in which efficiency
  is important.
  \begin{table}[htbp]
   \centering
  \caption{MSE of the LS,
ML, QML and EF estimators for
$\alpha_{01}=0.1,0.3$ with standard normal
errors}
\renewcommand{\arraystretch}{1.3}
\begin{tabular}{c|cc|cc|cc|cc}
\hline &\multicolumn{2}{|c|}{
$\hat{\alpha}_n^{(LS)}$}&\multicolumn{2}{|c|}{
$\hat{\alpha}_n^{(ML)}$}&\multicolumn{2}{|c}{
$\hat{\alpha}_n^{(QML)}$}
&\multicolumn{2}{|c}{$\hat{\alpha}_n^{(EF)}$}\\
\cline{2-9}
 Parameter&\multicolumn{2}{|c|}{
$n$}&\multicolumn{2}{|c|}{
$n$}&\multicolumn{2}{|c|}{
$n$}&\multicolumn{2}{|c}{ $n$}\\
 &50&500&50&500&50&500&50&500
 \\\hline
 &0.1195&0.0098& 0.0120
 &0.0013& 0.0120
 &0.0013&0.0062&0.0012 \\
 $\alpha_{01}=0.1$&0.7150& 0.0953&0.5668
 &0.0524&0.5668&0.0524
 &0.5670
 &0.0524\\
 &0.8345& 0.1051&0.5788&0.0537&0.5788&0.0537&0.5732
 &0.0536
\\\hline
 &0.1037&0.0978
 &0.0867&0.0420&0.0867&0.0420&0.0850
 &0.0405 \\
 $\alpha_{01}=0.3$&0.7554& 0.1899&0.6156
 &0.1324&0.6156&0.1324&0.5702
 &0.1320\\
 &0.8591& 0.2877&0.7023&0.1744&0.7023&0.1744&0.6552
 &0.1725
\\\hline
\end{tabular}
\newline{\small Note: The three
values in each cell are from top to bottom:
squared bias, variance, and MSE.}
\end{table}

\begin{table}[htbp]
   \centering
  \caption{MSE of the
LS, ML, QML and EF estimators for
$\alpha_{01}=0.1,0.3$ with $t_5$ distributed
errors}
\renewcommand{\arraystretch}{1.3}
\begin{tabular}{c|cc|cc|cc|cc}
\hline &\multicolumn{2}{|c|}{
$\hat{\alpha}_n^{(LS)}$}&\multicolumn{2}{|c|}{
$\hat{\alpha}_n^{(ML)}$}&\multicolumn{2}{|c}{
$\hat{\alpha}_n^{(QML)}$}
&\multicolumn{2}{|c}{$\hat{\alpha}_n^{(EF)}$}\\
\cline{2-9}
 Parameter&\multicolumn{2}{|c|}{
$n$}&\multicolumn{2}{|c|}{
$n$}&\multicolumn{2}{|c|}{
$n$}&\multicolumn{2}{|c}{ $n$}\\
 &50&500&50&500&50&500&50&500
 \\\hline
  &0.2314& 0.0124&0.0203&0.0088&
 0.0111&0.0050
 &0.0078
 &0.0013 \\
 $\alpha_{01}=0.1$&0.7808& 0.1090&0.6010&0.0552&
 0.5797&0.0537
 &0.5633
 &0.0556\\
 &1.0123& 0.1214&0.6213&0.0640&
 0.5908&0.0587
 &0.5711
 &0.0569
\\\hline
 &0.2132& 0.1096&0.1710&0.1013&0.0987
 &0.0268&0.0836
 &0.0277 \\
 $\alpha_{01}=0.3$&0.9875& 0.1711&0.6575&0.1579&0.6432
 &0.1410&0.6482
 &0.1377\\
 &1.2007& 0.2807&0.8285&0.2592&0.7419
 &0.1678&0.7318
 &0.1654
\\\hline
\end{tabular}
\newline{\small Note: The three
values in each cell are from top to bottom:
squared bias, variance, and MSE.}
\end{table}
    \newpage

\noindent\textbf{5. Proofs}\\ \\
 In this section we provide the proofs of Lemma 1, and Theorems 1 and
 2.\\\\
 \noindent {\sc Proof of Lemma 1}. Note from (3)
 that
  \begin{equation}\sqrt{n}(\hat{\theta}_n^{(LS)}-\theta_0)=
\left(\frac1n \sum_{t=1}^n Y_{t-1}
Y_{t-1}^T\right)^{-1} \left(n^{-1/2}\sum_{t=1}^n
Y_{t-1}\eta_t\right).
\end{equation}
 Since $\{X_t\}$ is stationary and ergodic, so is
 $(Y_{t-1}Y_{t-1}^T)$ which is finite by
 Assumption 1. Thus by the ergodic theorem
 \[
 \frac1n \sum_{t=1}^n Y_{t-1}
Y_{t-1}^T \stackrel{a.s.}{\longrightarrow}
\mathcal{U}.
 \]
\indent Next consider the second factor on the
right side of (14). Let $c=(c_0,\ldots,c_p)^T$ be
any vector with $c\not=0$. Recall that
$\eta_t=(\varepsilon_t^2-1)\sigma_t^2(\theta_0)$.
Then using the martingale central limit theorem
and Cramer-Wold device, it follows that
 \[
 n^{-1/2} \sum_{t=1}^n c^TY_{t-1}
\eta_t \stackrel{d}{\longrightarrow}
\mathcal{N}(0,Var(\epsilon_t^2)c^T\mathcal{R}(\theta_0)c).
 \]
  To prove this, note that $Y_{t-1}\eta_t$ is a
  martingale difference sequence since $E(Y_{t-1}
  \eta_t|\mathcal{F}_{t-1})=0$. We now verify the conditional Linderberg
  condition only since the other conditions can
  be verified easily.\\
  \indent For all $\epsilon>0$, we show that
\begin{eqnarray}
L_n&= &\frac1n \sum_{t=1}^n
(\sigma_t^2(\theta_0)c^TY_{t-1}
)^2\nonumber\\
&&\times E\{ (\varepsilon_t^2-1)^2
I(|(\varepsilon_t^2-1)\sigma_t^2(\theta_0)c^TY_{t-1}|>
\sqrt{n}\epsilon^2)|\mathcal{F}_{t-1}\}\nonumber\\
&=&o_p(1),
\end{eqnarray}
where $I(\Omega)$ is the indicator function of
the event $\Omega$. Observe that
 \begin{eqnarray*}
L_n&\le&\frac1n \sum_{t=1}^n
(\sigma_t^2(\theta_0)c^TY_{t-1}
)^2\nonumber\\
&&\times\, E\{ (\varepsilon_t^2-1)^2[
I(|\varepsilon_t^2-1|> n^{1/4} \epsilon)+I(
|\sigma_t^2(\theta_0)c^TY_{t-1}|>n^{1/4} \epsilon
) ]|\mathcal{F}_{t-1}\}\nonumber\\
&=&\frac1n \sum_{t=1}^n
(\sigma_t^2(\theta_0)c^TY_{t-1} )^2 E\{
(\varepsilon_t^2-1)^2
I(|\varepsilon_t^2-1|> n^{1/4} \epsilon)\}\nonumber\\
&&+\,
 Var(\varepsilon_t^2)\frac1n \sum_{t=1}^n (\sigma_t^2(\theta_0)c^TY_{t-1}
)^2
 I( |\sigma_t^2(\theta_0)c^TY_{t-1}|>n^{1/4}
\epsilon )\nonumber\\
 &\equiv& T_1+  Var(\varepsilon_t^2) T_2.
\end{eqnarray*}
  Noting that $Var(\varepsilon_t^2)<\infty$ and
  that
  \[
  \frac1n \sum_{t=1}^n (\sigma_t^2(\theta_0)c^TY_{t-1}
)^2=E[(\sigma_t^2(\theta_0)c^TY_{t-1}
)^2]+o_p(1),
   \]
   we obtain $T_1=o_p(1)$. Moreover, note from
   Assumption 1 that
   \[
   E(T_2)=E[(\sigma_t^2(\theta_0)c^TY_{t-1}
)^2I( |\sigma_t^2(\theta_0)c^TY_{t-1}|>n^{1/4}
\epsilon )]=o(1),
   \]
   which implies $T_2\ge 0$. Hence (15) is
   satisfied and by Slutsky's theorem the
   assertion of Lemma 1 follows.\\\\
\noindent {\sc Proof of Theorem 1}. From (8),
observe that
\[\sqrt{n}(\hat{\theta}_n^{(EF)}-\theta_0)=
\left(\frac1n \sum_{t=1}^n
\frac{Y_{t-1}
Y_{t-1}^T}{((\hat{\theta}_n^{(LS)})^T
Y_{t-1})^2}\right)^{-1}
\left(n^{-1/2}\sum_{t=1}^n
\frac{Y_{t-1}\eta_t}{((\hat{\theta}_n^{(LS)})^T
Y_{t-1})^2}\right).
\]
The result in the theorem can be proved if we
show that
\begin{equation}
\frac1n \sum_{t=1}^n Y_{t-1} Y_{t-1}^T
\{((\hat{\theta}_n^{(LS)})^T Y_{t-1}Y_{t-1}^T
\hat{\theta}_n^{(LS)})^{-1}-(\theta_0^T
Y_{t-1}Y_{t-1}^T \theta_0)^{-1}\}=o_p(1)
\end{equation}
and
\begin{equation}
n^{-1/2} \sum_{t=1}^n Y_{t-1} \eta_t
\{((\hat{\theta}_n^{(LS)})^T Y_{t-1}Y_{t-1}^T
\hat{\theta}_n^{(LS)})^{-1}-(\theta_0^T
Y_{t-1}Y_{t-1}^T \theta_0)^{-1}\}=o_p(1).
\end{equation}
In view of Lemma 1, we note that
$\hat{\theta}_n^{(LS)}\stackrel{a.s.}
{\rightarrow}\theta_0$ for sufficiently large
$n$, and thus,
$\sigma_t^2(\hat{\theta}_n^{(LS)})$ behaves like
$\sigma_t^2(\theta_0)$ for each $t=1,\ldots,n$.
This statement holds true, if for $\epsilon>0$,
there exists $N_\epsilon$ such that $\|
\hat{\theta}_n^{(LS)}-\theta_0\|\le \epsilon$ for
all $n> N_\epsilon$, with probability one.
Consequently, we have
\[\sqrt{n}(\hat{\theta}_n^{(EF)}-\theta_0)=
\left(\frac1n \sum_{t=1}^n
\frac{Y_{t-1} Y_{t-1}^T}{(\theta_0^T
Y_{t-1})^2}\right)^{-1}
\left(n^{-1/2}\sum_{t=1}^n
\frac{Y_{t-1}\eta_t}{(\theta_0^T
Y_{t-1})^2}\right)+o_p(1),
\]
for which, the proof is reduced to that of Lemma
1. Since the proof of (16) and (17) is similar,
we
prove (16) only as follows. \\
\indent By a Taylor expansion around
$\hat{\theta}_n^{(LS)}$ at $\theta_0$, it follows
that (16) is dominated by
\[
O_p(\sqrt{n}(\hat{\theta}_n^{(LS)}-\theta_0))\times
\biggl| \frac1n \sum_{t=1}^n Y_{t-1} \eta_t
\tilde{\theta}_\epsilon^T
(Y_{t-1}Y_{t-1}^T)(\tilde{\theta}_\epsilon^T
Y_{t-1}Y_{t-1}^T
\tilde{\theta}_\epsilon)^{-2}\biggl|,
\]
where $\tilde{\theta}_\epsilon$ lies between
$\theta_0$ and $\hat{\theta}_n^{(LS)}$. Moreover,
from the ergodic theorem, we readily see that
\[
\frac1n \sum_{t=1}^n Y_{t-1} \eta_t
\tilde{\theta}_\epsilon^T
(Y_{t-1}Y_{t-1}^T)(\tilde{\theta}_\epsilon^T
Y_{t-1}Y_{t-1}^T
\tilde{\theta}_\epsilon)^{-2}\stackrel{a.s.}
{\rightarrow} 0.
\]
Hence the conception of Theorem 1 follows.\\\\
\noindent{\sc Proof of Theorem 2.} The proof of
this theorem requires the following matrix
inequality (see Kholevo (1969)).\\\\
\textbf{Lemma 3.} {\it Let $A(u)$ and $B(u)$ be
$r\times s$ and $t\times s$ random matrices,
respectively, and $\psi(u)$ is a function that is
positive everywhere. If $E(BB^T/\psi)$ exists,
then
\[
E(AA^T \psi)\ge E(A B^T) E(BB^T/\psi)^{-1}E[(A
B^T )]^T.
\]
The equality holds if and only if there exists a
constant $r\times t$ matrices $C$ such that $\psi
A+CB=0$ almost everywhere.}
\\ \\
\indent Using the notation of Lemma 3, write
$A=B=Y_{t-1}$ and $\psi=\sigma_t^4(\theta_0)$.
Then it follows that $\mathcal{R}(\theta_0) \ge
\mathcal{U} \Gamma^{-1}(\theta_0) \mathcal{U}$.
It is obvious from Lemma 3 that the equality
$\mathcal{R}(\theta_0)= \mathcal{U}
\Gamma^{-1}(\theta_0) \mathcal{U}$ holds if and
only if $\sigma_t^2(\theta_0)=k$, a constant
almost everywhere. Hence we get the desired
result.\\\\
\textbf{6. References}
\begin{itemize}
\item[]  Bauwens, L., Laurent, S. and Rombouts,
J. (2006). Multivriate GARCH Models: A survey.
{\it J. Appl. Econ.} 21, 79-109.
 \item []
 Bollerslev, T. (1986). Generalized autoregressive
 conditional heteroskedasticity. {\it J.
 Econometrics} 31, 307-328.
    \item [] Bougerol, P. and Picard, N. (1992a).
Strict stationarity of generalized autoregressive
processes. {\it Ann. Prob.} 20, 1714-1730.
 \item []Bougerol, P. and Picard, N. (1992b).
Stationarity of GARCH processes and of
 some nonnegative time series. {\it J. Econometrics}
 52, 115-127.

    \item [] Engle, R. F. (1982). Autoregressive conditional
 heteroskedasticity with estimates of the variance
 of U.K. inflation. {\it Econometrica} 50, 987-1007.
    \item [] Engle, R. (1995). {\it ARCH Selected
Readings}. New York: Oxford University Press.
    \item [] Francq, C. and Zako\"{\i}an, J. M.
    (2004). Maximum likelihood estimation of pure
    GARCH and ARMA-GARCH processes. {\it
    Bernoulli} 10(4), 605-637.
    \item [] Godambe, V. P. (1985). The foundations of finite sample
 estimation in stochastic processes. {\it Biometrika}, {\bf 72}, 419-28.
    \item [] Gouri\'eroux, C. (1997). {\it
ARCH Models and Financial Applications}. New
York: Springer. \item [] Hansen L. P. (1982).
Large sample properties of generalised method of
moments estimators. {\it Econometrica} {\bf 50},
1029-1054. \item [] Hosoya, Y. and Tanguchi, M.
(1982). A central limit theorem for stationary
processes and the parameter estimation of linear
processes. {\it Ann. Statist.} {\bf 10} 132-153.
Correction: (1993). {\bf 21} 1115-1117. \item []
Kholevo, A.S. (1960). On estimates of regression
coefficients. {\it Theory Prob. Appl.} 14,
79-104. \item [] Li, D. X. and Turtle, H. J.
(2000). Semiparametric ARCH models: an estimating
function approach. {\it J. Bus. Econ. Statist.}
18(2), 174-186. \item [] Mikosch, T. (2003).
Modeling dependence and tails of financial time
series. In {\it Extreme Values in Finance,
Telecommunications, and the Environment} (B.
Finkenst\"adt and H. Rootz\'en, eds.), 185-286,
Chapman and Hall, Boca Raton, FL.
    \item [] Nelson, D. B. (1990). Stationarity and persistence in the
 GARCH(1,1) model. {\it Econ. Theory} 6, 318-34.
\item [] Rich, R. W., Raymond, J. and Butler, J.
S. (1991). Generalized instrumental variables
estimation of autoregressive conditional
heteroskedastic models. {\it Economic Letters}
35, 179-185.
\item [] Shephard, N. (1996).
Statistical aspects of ARCH and stochastic
volatility. In D. R. Cox, D. V. Hinkley and O. E.
Barndorff-Nielson (eds), {\it Time Series Models
in Econometrics, Finance and Other Fields}, 1-55.
London: Chapman $\&$ Hall.
    \item [] Tj$\o$stheim, D. (1986). Estimation in nonlinear
time
 series models. {\it Stoch. Proc. Appl.}, {\bf 21}, 251-273.
    \item [] Weiss, A. A. (1986). Asymptotic
    theory for ARCH models: Estimation and
    testing. {\it Econometric Theory} 2, 107-31.
\end{itemize}

\end{document}